\documentclass[12pt]{amsart}
\usepackage{amsmath,amscd,amssymb,amsfonts}
\newcommand{\ver}{June 2, 2006, v.2}
\newcommand{\scirc}{\,\raise.2ex\hbox{${\scriptstyle\circ}$}\,}
\newcommand{\sbull}{\,\raise.2ex\hbox{${\scriptstyle\bullet}$}\,}
\newcommand{\mprod}{\hbox{$\prod$}}
\newcommand{\msum}{\hbox{$\sum$}}
\newcommand{\mopl}{\hbox{$\bigoplus$}}
\newcommand{\mtim}{\hbox{$\times$}}
\newcommand{\motim}{\hbox{$\otimes$}}
\newcommand{\bQ}{{\mathbb Q}}
\newcommand{\bR}{{\mathbf R}}
\newcommand{\bZ}{{\mathbb Z}}
\newcommand{\cF}{{\mathcal F}}
\newcommand{\cH}{{\mathcal H}}
\newcommand{\cO}{{\mathcal O}}
\newcommand{\cS}{{\mathcal S}}
\newcommand{\cZ}{{\mathcal Z}}
\newcommand{\ok}{\overline{k}}
\newcommand{\oK}{\overline{K}}

\newcommand{\tf}{\widetilde{f}}
\newcommand{\tS}{\widetilde{S}}

\newcommand{\tpi}{\widetilde{\pi}}
\newcommand{\charr}{\hbox{\rm char}\,}
\newcommand{\Alb}{\hbox{\rm Alb}}
\newcommand{\CH}{\hbox{\rm CH}}
\newcommand{\IC}{\hbox{\rm IC}}
\newcommand{\End}{\hbox{\rm End}}
\newcommand{\Ext}{\hbox{\rm Ext}}
\newcommand{\Hom}{\hbox{\rm Hom}}
\newcommand{\Ker}{\hbox{\rm Ker}}
\newcommand{\rIm}{\hbox{\rm Im}}
\newcommand{\Spec}{\hbox{\rm Spec}}
\newcommand{\cSpec}{{\mathcal S}pec}
\newcommand{\Cor}{\hbox{\rm Cor}}
\newcommand{\red}{\text{\rm red}}

\begin{document}
\title[Chow-K\"unneth decomposition]
{Chow-K\"unneth decomposition for varieties\\
with low cohomological level}
\author{Morihiko Saito}
\address{RIMS Kyoto University, Kyoto 606-8502 Japan}
\email{msaito@kurims.kyoto-u.ac.jp}
\date{\ver}
\begin{abstract}
We show that a smooth projective variety admits a Chow-K\"unneth
decomposition if the cohomology has level at most one except for
the middle degree.
This can be extended to the relative case in a weak sense if the
morphism has only isolated singularities, the base space is
1-dimensional, and the generic fiber satisfies the above condition.
\end{abstract}
\maketitle

\centerline{\bf Introduction}

\bigskip\noindent
Let $X$ be a smooth projective variety over a perfect field $k$.
Murre's conjecture ([16], [17]) predicts the existence of projectors
$\pi_j$ in the Chow group of $X\mtim_kX$ with rational coefficients
which are orthogonal to each other and whose images in the Chow group
modulo homological equivalence give the usual K\"unneth decomposition
of the diagonal class.
So far this conjecture is proved for curves, surfaces, abelian
varieties, uniruled 3-folds, Calabi-Yau 3-folds, etc. (see [12]).

Let $l$ be a prime different from the characteristic of $k$.
Let $X_{\ok}$ be the base change of $X$ by an algebraic closure
$\ok$ of $k$.
We say that $H^j(X_{\ok},\bQ_l)$ has level $\le i$ if there is
a closed subvariety $Z$ of codimension $p\ge (j-i)/2$ such that
the restriction morphism $H^j(X_{\ok},\bQ_l)\to H^j(X_{\ok}\setminus
Z_{\ok},\bQ_l)$ vanishes.
In this paper we prove

\medskip\noindent
{\bf Theorem~1.}
{\it Assume $H^j(X_{\ok},\bQ_l)$ has level $\le 1$ for any
$j<n:=\dim X$.
Then $X$ admits a Chow-K\"unneth decomposition.
}

\medskip\noindent
{\bf Corollary.}
{\it If $\dim X = 3 $, $\Gamma(X,\Omega_X^2)=0$ and $\charr k = 0$,
then $X$ admits a Chow-K\"unneth decomposition.
}

\medskip
For the proof of Theorem~1 we introduce {\it projectors of
restricted type} to make use of Murre's argument [16] in a most
effective way.
This simplifies especially the arguments about the orthogonality
of the projectors.
For a generalization to the higher level case, we would need
stronger conjectures of Murre [17] for certain varieties of
lower dimensions together with the standard conjectures [13].
In case of general threefolds, we would need the standard
conjecture of Lefschetz-type for $X$ together with the
injectivity of $\End(h^2(Y))\to\End(H^2(Y_{\ok},\bQ_l))$,
where $Y$ is a general hyperplane section of $X$ and
$h^j(Y)$ is the `image' of the projector $\pi_j$.

If the condition on level is satisfied also for $j=n$,
one might define the projector $\pi_n$ in the same way
as in the proof of Theorem~1.
However, it is quite nontrivial whether we have the equality
$\sum_j\pi_j=\Delta_X$ in $\CH^n(X\mtim_kX)_{\bQ}$
(where $\Delta_X$ is the diagonal of $X$), because the difference
may correspond to a `fantom' motive, see [12], [18].
This is closely related to Bloch's conjecture [3] for surfaces with
$p_g=0$.

Using the theory of relative correspondences [4],
we can generalize Theorem~1 to the relative case in a weak sense
as follows:

\medskip\noindent
{\bf Theorem~2.} {\it
Let $f:X\to S$ be a surjective morphism of irreducible smooth
projective varieties over a field $k$ of characteristic $0$
with $\dim S=1$.
Let $Y_K$ be the generic fiber of $f$ over $K=k(S)$.
Assume $f$ has only isolated singularities.
If $H^j(Y_{\oK},\bQ_l)$ has level $\le 1$ for $j<m:=\dim Y_K$,
then $f$ admits a relative Chow-K\"unneth decomposition
in the weak sense {\rm(}see {\rm(3.3))}.
If $H^j(Y_{\oK},\bQ_l)$ has level $\le 1$ for $j\le m$, then
$X$ admits a Chow-K\"unneth decomposition.
}

\medskip
Similar results were obtained in [9], [10] where the assumptions
on the generic and special fibers are quite different from ours.
It is not easy to prove the relative Chow-K\"unneth decomposition
in the strong sense (see [9]) in our situation unless $m=1$,
although it is easy to construct the sum of the projectors
corresponding the direct factors of the direct image with
discrete supports.

In Section 1 we introduce projectors of restricted type.
In Section 2 we prove Theorem~1 using Murre's construction
together with projectors of restricted type.
In Section 3 we review the theories of relative correspondences
and relative Chow-K\"unneth decompositions.
In Section 4 we study the case of relative dimension 1.
In Section 5 we prove Theorem~2.

\bigskip\bigskip
\noindent\centerline{\bf 1. Projectors of restricted type}

\bigskip\noindent
{\bf 1.1.~Definition.}
Let $X$ be a smooth irreducible projective variety over a perfect
field $k$, and $n=\dim X$.
We say that a projector $\pi\in\CH^n(X\mtim_kX)_{\bQ}$ has pure
cohomological degree $j$ if its action on $H^{i}(X_{\ok},\bQ_l)$
vanishes for $i\ne j$.
(Here a projector means $\pi^2=\pi$ and the composition of
correspondences are defined as usual.)
We say that $\pi$ is a {\it projector of restricted type} with 
degree $j$ if there is a smooth projective $k$-variety $Y$ of
pure dimension $m$ having a projector $\pi'$ with pure
cohomological degree $j'=j+2i$ together with correspondences
$\Gamma\in\CH^{n+i}(X\mtim_kY)_{\bQ}$,
$\Gamma'\in\CH^{m-i}(Y\mtim_kX)_{\bQ}$ such that
$$\aligned &\pi=\Gamma'\scirc\pi'\scirc\Gamma,\\
&\rIm(\Gamma_*:H^{j}(X_{\ok},\bQ_l)(-i)\to H^{j'}(Y_{\ok},\bQ_l))
=\rIm\,\pi'_*,\\&\Ker(\Gamma'_*:H^{j'}(Y_{\ok},\bQ_l)\to H^{j}
(X_{\ok},\bQ_l)(-i))=\Ker\,\pi'_*,\endaligned\leqno(1.1.1)$$
and furthermore we have the inclusion
$$\End(\rIm\,\pi')\hookrightarrow\End(H^{j'}(X_{\ok},\bQ_l)).
\leqno(1.1.2)$$

\medskip\noindent
{\bf 1.2.~Remarks.} (i)
By the definition of the `image' of a projector, $\End(\rIm\,\pi'))$
consists of the correspondences of the form
$\pi'\scirc\xi\scirc\pi'$.

\medskip
(ii)
If $j'=1$, the inclusion (1.1.2) in the above definition follows
from the theory of abelian varieties ([14], [15], [20]), see [16], [19].
In general, it is closely related to Murre's conjectures, see [17].

\medskip
(iii)
The above conditions imply
$$\aligned H^{j'}(Y_{\ok},\bQ_l)&=\Ker\,\Gamma'_*\oplus\rIm\,\Gamma_*
\\H^{j}(X_{\ok},\bQ_l)&=\Ker\,\Gamma_*\oplus\rIm\,\Gamma'_*,
\endaligned$$ together with $\Ker\,\Gamma_*=\Ker\,\pi_*$,
$\rIm\,\Gamma'_*=\rIm\,\pi_*\simeq\rIm\,\pi'_*$, and
$$\Gamma_*\scirc\Gamma'_*=\pi'_*\quad\text{on}\quad
H^{j'}(Y_{\ok},\bQ_{l}).\leqno(1.2.1)$$
Then we have by (1.1.2)
$$\pi'\scirc\Gamma\scirc\Gamma'\scirc\pi'=\pi',\leqno(1.2.2)$$
although it is not clear whether $\Gamma\scirc\Gamma'=\pi'$ in
$\CH^m(Y\mtim_kY)_{\bQ}$.

\medskip\noindent
{\bf 1.3.~Lemma.}
{\it Let $\xi,\xi'\in\CH^n(X\mtim_kX)_{\bQ}$ whose actions on
$H^j(X_{\ok},\bQ_{l})$ are the identity.
Let $\pi$ be a projector of restricted type with degree $j$.
Then so is $\tpi:=\xi'\scirc\pi\scirc\xi$.
}

\medskip\noindent
{\it Proof.}
We have $\tpi^2=\tpi$ by (1.1.2) together with (1.2.1).

\medskip\noindent
{\bf 1.4.~Proposition.}
{\it Let $\pi$ be a projector of restricted type with degree $j$.
Let $\pi'_i$ be projectors which commute with each other and
whose actions on $H^j(X_{\ok},\bQ_{l})$ are zero for $i\in I$.
Set $\tpi=\prod_{i\in I}(1-\pi'_i)\scirc\pi\scirc\prod_{i\in I}
(1-\pi'_i)$. Then $\pi'$ is a projector of restricted type
with degree $j$ which is orthogonal to the $\pi'_i$, and
$\rIm\,\tpi_*=\rIm\,\pi_*$ in $H^j(X_{\ok},\bQ_{l})$.
}

\medskip\noindent
{\it Proof.}
This follows from Lemma~(1.3).

\bigskip\bigskip
\noindent\centerline{\bf 2. Construction of projectors}

\bigskip\noindent
{\bf 2.1.~Proposition.} {\it
Let $X$ be a smooth projective variety over a perfect filed $k$.
Fix an integer $j$,
and assume that the level of $H^j(X_{\ok},\bQ_l)$ is at most $j-2p$.
Then there is a smooth projective variety $Y$ of pure dimension
$n-p$ together with a morphism $\rho:Y\to X$ inducing the surjection
$$(^t\Gamma_{\rho})_*:H^{j-2p}(Y_{\ok},\bQ_l)(-p)\to H^j(X_{\ok}
\bQ_l),\leqno(2.1.1)$$ where $\Gamma_{\rho}$ is the algebraic cycle
defined by the graph of $\rho$, and $(^t\Gamma_{\rho})_*$ coincides
with the Gysin morphism $\rho_*$.
}

\medskip\noindent
{\it Proof.}
By assumption there is a closed subvariety $Z$ of codimension
$p$ such that we have the vanishing of the restriction morphism
$$H^j(X_{\ok},\bQ_l)\to H^j(X_{\ok}\setminus Z_{\ok},\bQ_l).$$
Replacing $Z$ if necessary, we may assume that $Z$ is pure
dimensional and furthermore it is a complete intersection.
Let $Y$ be the disjoint union of smooth projective alterations of
the irreducible components of $Z$, see [5].
Then, using the above condition on $Z$ we get the surjectivity of
(2.1.1).
(Indeed, this can be proved by using the decomposition theorem [2]
because $\bR\Gamma_Z\bQ_l$ is a perverse sheaf up to a shift
and its lowest weight part is the direct sum of the intersection
complexes of the irreducible components of $Z$.)

\medskip\noindent
{\bf 2.2.~Albanese and Picard varieties.}
Let $X,Y$ be smooth projective varieties over a field $k$, and
$\Alb_X,P_Y$ denote respectively the Albanese and Picard varieties
of $X,Y$. Then we have a canonical isomorphism (see e.g. [16], [19])
$$\Hom(\Alb_X,P_Y)\motim\bQ=\CH^1(X\mtim_kY)_{\bQ}/
(pr_1^*\CH^1(X)_{\bQ}+pr_2^*\CH^1(Y)_{\bQ}).$$
Note that this is identified with a subspace of $H^1(X_{\ok},\bQ_l)
\motim H^1(Y_{\ok},\bQ_l)$.

\medskip\noindent
{\bf 2.3.~Proposition.}
{\it With the notation and the assumption of {\rm (2.1)}, assume
furthermore $j=2p+1$.
Then we can replace $Y$ and the cycle ${}^t\Gamma_{\rho}$ with a
purely 1-dimensional smooth projective variety $C$ and
$\Gamma\in\CH^{p+1}(C\mtim_kX)_{\bQ}$.
}

\medskip\noindent
{\it Proof.}
Let $C$ be the intersection of general hyperplane sections of $Y$.
More precisely, it is the fiber of $\cZ\to\cS$ over a general
closed point of $\cS$, where $\cS$ is the $r$-ple product of
the parameter space of hyperplanes of $Y$ (where $r=\dim Y-1$),
and $\cZ$ is a closed subvariety of $Y\mtim_k\cS$ whose fiber
over a geometric point of $\cS$ is the intersection of the
corresponding $r$ hyperplane sections.

By (2.2) together with the hard Lefschetz theorem [8],
there is $\Lambda\in\CH^1(Y\mtim_kY)_{\bQ}$ such that the
action of $\Lambda\scirc{}^t\Gamma_{\iota}\scirc\Gamma_{\iota}$ on
$H^1(Y_{\ok},\bQ_l)$ is the identity, where $\Gamma_{\iota}$ is the
graph of the inclusion $\iota:C\to Y$, see [16], [19].
Then $\Gamma$ is defined by
${}^t\Gamma_{\rho}\scirc\Lambda\scirc{}^t\Gamma_{\iota}$.

\medskip\noindent
{\bf 2.4.~Proof of Theorem~1.}
We construct the projectors $\pi_j$ and $\pi_{2n-j}$ by increasing
induction on $j<n$.
For $j=n$, put $\pi_n=\Delta_X-\sum_{j\ne n}\pi_j$ as usual.
Since the assertion is proved by Murre [16] for $j\le 1$,
we may assume $j\ge 2$.
Here we have either $j=2p$ or $2p+1$.

We fist consider the case $j=2p+1$.
By Proposition (2.3) there is a curve $C$ together with
$\Gamma\in\CH^{p+1}(C\mtim_kX)_{\bQ}$ inducing the surjection
$$\Gamma_*:H^1(C_{\ok},\bQ_l)\to H^j(X_{\ok},\bQ_l)(p).
\leqno(2.4.1)$$
Let $\Gamma_H\in\CH^{n+1}(X\mtim_kX)$ be the correspondence defined
by the diagonal of a general hyperplane section of $X$. Consider
$$(\Gamma_H^{n-j}\scirc\Gamma)_*:H^1(C_{\ok},\bQ_l)\to
H^{2n-j}(X_{\ok},\bQ_l)(n-p-1).$$
This is the composition of (2.4.1) with $L^{n-j}$, and is hence
surjective by the hard Lefschetz theorem [8],
where $L$ is the cup product with the hyperplane section class.
Taking the dual, we get the injectivity of
$$({}^t\Gamma\scirc\Gamma_H^{n-j})_*:H^j(X_{\ok},\bQ_l)(p)\to
H^1(C_{\ok},\bQ_l),\leqno(2.4.2)$$
because $\Gamma_H$ is self-dual.
Consider the composition of (2.4.1) and (2.4.2)
$$({}^t\Gamma\scirc\Gamma_H^{n-j}\scirc\Gamma)_*:
H^1(C_{\ok},\bQ_l)\to H^1(C_{\ok},\bQ_l).$$

The correspondence ${}^t\Gamma\scirc\Gamma_H^{n-j}\scirc\Gamma\in
\CH^1(C\mtim_kC)_{\bQ}$ induces a morphism of abelian varieties
$\phi:P_{C}\to\Alb_C$ in the notation of (2.2).
(Here $P_{C}=\Alb_C$ because $C$ is a curve.)
By the semisimplicity of abelian varieties there is a morphism
$\psi:\Alb_C\to P_C$ such that $\psi\scirc\phi$ is an idempotent
up to a nonzero constant multiple, and
$$\Ker\,\psi\scirc\phi=\Ker\,\phi,\quad\rIm\,\psi\scirc\phi=
\rIm\,\psi.\leqno(2.4.3)$$
Indeed, this is obtained by splitting the short exact
sequences associated to $\Ker\,\phi$, $\rIm\,\phi$, etc.
(up to an isogeny).

By (2.2), $\psi$ corresponds to a correspondence
$\Lambda\in\CH^1(C\mtim_kC)_{\bQ}$ whose cohomology class is
contained in $H^1(C_{\ok},\bQ_l)\motim H^1(C_{\ok},\bQ_l)(1)$.
Set
$$\pi'_C:=\Lambda \scirc{}^t\Gamma\scirc\Gamma_H^{n-j}\scirc
\Gamma\in\CH^1(C\mtim_kC)_{\bQ}.$$
Replacing $\Lambda$ with $\alpha\Lambda$ for an appropriate
$\alpha\in\bQ^*$, $\pi'_C$ is a projector with pure cohomological
degree $1$ because $\Lambda \scirc{}^t\Gamma\scirc\Gamma_H^{n-j}
\scirc\Gamma\scirc\Lambda=\Lambda$ by (2.2).
Moreover we have by (2.4.2--3)
$$\Ker\,(\pi'_C)_*=\Ker\,\Gamma_*,\,\,\,\rIm\,(\pi'_C)_*=
\rIm\,(\Lambda\scirc{}^t\Gamma\scirc\Gamma_H^{n-j})_*\,\,\,
\text{in}\,\,\, H^1(C_{\ok},\bQ_l).$$
Let
$$\pi'_j=\Gamma\scirc\pi'_C\scirc\Lambda\scirc{}^t\Gamma
\scirc\Gamma_H^{n-j}.$$
Then this is a projector of restricted type with degree $j$.
Set
$$\pi''_j=\mprod_{|n-i|>n-j}(1-\pi_i)\scirc\pi'_j\scirc
\mprod_{|n-i|>n-j}(1-\pi_i).$$
By Proposition (1.4) this is also a projector of restricted type with
degree $j$, and is orthogonal to $\pi_i$, $\pi_{2n-i}$ for $i<j$.
In order to define $\pi_{2n-j}$ by ${}^t\pi_j$, we have to modify
$\pi''_j$ so that it becomes orthogonal to its transpose.
Here we have ${}^t\pi''_j\scirc\pi''_j=0$, because
${}^t\pi''_j\scirc\pi''_j$ factors through the composition of
correspondences
${}^t\Gamma\scirc u\scirc\Gamma\in\CH^{1+j-n}(C\mtim_kC)_{\bQ}$
whose cohomology class vanishes in the case $n=j+1$, where $u$ is
the product of $(1-{}^t\pi_i)(1-\pi_i)$ for $|n-i|>n-j$.
Then we can replace $\pi''_j$ with $\pi_j=\pi''_j\scirc
(1-{}^t\pi''_j/2)$ so that we get $\pi_j\scirc{}^t\pi_j=0$ as in [19].

In the case $j=2p$, the curve $C$ is replaced by a $0$-dimensional
variety, and abelian varieties are replaced by finite dimensional
$\bQ$-vector spaces. Then the argument is similar (and easier).
This finishes the proof of Theorem~1.

\medskip\noindent
{\bf 2.5.~Remark.}
Let $j$ be an integer in $[0,n-1]$ where $n=\dim X$.
Let $Y$ be an intersection of general hyperplane sections with
$\dim Y=j$.
Let $\iota:Y\to X$ denote the inclusion morphism. 
If the standard conjecture of Lefschetz-type [13] holds for $X$,
there is a correspondence $\Lambda\in\CH^j(X\mtim_kX)_{\bQ}$
such that the action of $\Lambda\scirc{}^t\Gamma_{\iota}\scirc
\Gamma_{\iota}$ on $H^j(X_{\ok},\bQ_l)$ is the identity, where
$\Gamma_{\iota}\in\CH^n(X\mtim_kY)$ is defined by the graph of
$\iota$.
If furthermore we have the injectivity 
$$ \End(h^j(Y))\hookrightarrow\End(H^j(Y_{\ok},\bQ_l)),
\leqno(2.5.1)$$ where $h^j(Y)$ is the image of the projector
$\pi_j^Y$ associated to some Chow-K\"unneth decomposition for
$Y$, then $\pi_j^Y\scirc\Gamma_{\iota}\scirc\Lambda\scirc{}^t
\Gamma_{\iota}\scirc\pi_j^Y$ is an idempotent, and
the projector $\pi_j$ for $X$ can be defined by
modifying $\Lambda\scirc{}^t\Gamma_{\iota}\scirc\pi_j^Y\scirc
\Gamma_{\iota}$ as in (1.4).

These two hypotheses are satisfied for $j=0,1$.
We can choose any order of $\{0,\dots,n-1\}$ to proceed by
induction on $j$ as long as the hypotheses are satisfied for
each $j$.
At the first step of the induction, the second assumption can
be replaced by a weaker condition that
$\pi_j^Y\scirc\Gamma_{\iota}\scirc\Lambda\scirc{}^t
\Gamma_{\iota}\scirc\pi_j^Y$ is an idempotent,
because we do not need the argument for the orthogonality with
the other projectors.
However, this is still quite difficult in the case of threefolds
with $j=2$.

\bigskip\bigskip
\centerline{\bf 3. Relative decomposition}

\bigskip\noindent
{\bf 3.1.~Relative correspondences.}
Let $X,Y$ be a smooth $k$-varieties with proper morphisms
$f:X\to S,g:Y\to S$.
The group of relative correspondences of degree
$j$ between $X$ and $Y$ over $S$ (see [4]) is defined by
$$\Cor_S^j(X,Y)=\CH_{\dim Y-j}(X\mtim_SY)_{\bQ},$$
where $Y$ is assumed to be pure-dimensional.
In general, we define $\Cor_S^j(X,Y)$ to be the direct sum of
$\Cor_S^j(X,Y_i)$ where $Y_i$ are the irreducible components of
$Y$.
The composition
$$\Cor_S^j(X,Y)\times\Cor_S^i(Y,Z)\to\Cor_S^{i+j}(X,Z)
\leqno(3.1.1)$$
is defined by using the pull-back by the morphism
$$X\mtim_SY\mtim_SZ\to(X\mtim_SY)\mtim_k(Y\mtim_SZ),$$
and then the push-forward by the morphism to $X\mtim_SZ$.
Note that these morphisms are respectively the base change
of the diagonal morphism $Y\to Y\mtim_kY$ and that of $Y\to S$.
The composition of $\xi\in\Cor_S^j(X,Y)$ and $\eta\in\Cor_S^i
(Y,Z)$ is usually denoted by $\eta\scirc\xi$.
We will denote it also by $\xi\sbull\eta$.

By a similar argument we have the action of the correspondences
$$\Cor_S^j(X,Y)\times\CH^i(X)_{\bQ}\to\CH^{i+j}(Y)_{\bQ}.
\leqno(3.1.2)$$
This is essentially a special case of (3.1.1) where $X,Y,Z$ are
replaced by $S,X,Y$ respectively.

By [4] each $\gamma\in\Cor_S^j(X,Y)$ induces natural morphisms
$$\aligned
\gamma_*:\bR(f_{\ok})_*\bQ_l&\to\bR(g_{\ok})_*\bQ_l(j)[2j],\\
\gamma_*:\bR\Gamma(X_{\ok},\bQ_l)&\to\bR\Gamma(Y_{\ok},
\bQ_l)(j)[2j],\endaligned\leqno(3.1.3)$$
in a compatible way with compositions,
where $f_{\ok}:X_{\ok}\to S_{\ok},g_{\ok}:Y_{\ok}\to S_{\ok}$
are the base changes of $f,g$.
We have a natural morphism
$$\Cor_S^j(X,Y)\to\Cor_k^j(X,Y),\leqno(3.1.4)$$
and a relative correspondence naturally induces an absolute
correspondence.
This is compatible with compositions and also with (3.1.3).

\medskip\noindent
{\bf 3.2.~Decomposition theorem.}
Let $f:X\to S$ be a proper morphism of irreducible varieties
over $k$ such that $X$ is smooth over $k$. Set $n=\dim X$
so that $\bQ_l[n]$ is a perverse sheaf on $X_{\ok}$.
Let $f_{\ok}:X_{\ok}\to S_{\ok}$ denote the base change of $f$
by $k\to\ok$.
By [2] there is a non-canonical isomorphism
$$\bR(f_{\ok})_*\bQ_l[n]\simeq\mopl_j({}^pR^j(f_{\ok})_*
\bQ_l[n])[-j]\quad\text{in}\,\,\, D_c^b(S_{\ok},\bQ_l),
\leqno(3.2.1)$$
together with canonical isomorphisms in the category
of perverse sheaves
$${}^pR^j(f_{\ok})_*\bQ_l[n]=\mopl_{Z}\IC_ZE^j_{Z^{\circ}},
\leqno(3.2.2)$$
where $Z$ runs over integral closed subvarieties of
$S$, and $E^j_{Z^{\circ}}$ is a smooth $\bQ_l$-module
on $Z^{\circ}_{\ok}$ with $Z^{\circ}$ a dense open subvariety
of $Z$ such that $E^j_{Z^{\circ}}=0$ except for a finite
number of $Z$.
Here $\IC_ZE^j_{Z^{\circ}}$ denotes the intersection complex,
i.e. the intermediate direct image of $E_{Z^{\circ}}[\dim Z]$.
Note that ${}^pR^j(f_{\ok})_*\bQ_l[n]={}^p\cH^j\bR(f_{\ok})_*
\bQ_l[n]$ with ${}^p\cH^{j}$ the perverse cohomology functor,
where ${}^pR^j(f_{\ok})_*\bQ_l[n]$ means ${}^pR^j(f_{\ok})_*
(\bQ_l[n])$ and not $({}^pR^j(f_{\ok})_*\bQ_l)[n]$.
By (3.2.1--2) we get a non-canonical isomorphism
$$\bR(f_{\ok})_*\bQ_l[n]\simeq\mopl_j\mopl_{Z}
(\IC_ZE^j_{Z^{\circ}})[-j].\leqno(3.2.3)$$
These assertions can be reduced to the case where $k$ is
finite using a model of $f$.

\medskip\noindent
{\bf 3.3.~Relative Chow-K\"unneth decomposition.}
With the notation and the assumption of (3.2),
we say that $f:X\to S$ admits a {\it relative Chow-K\"unneth
decomposition in the weak sense},
if there are mutually orthogonal idempotents
$\pi_i\in\Cor^0_S(X,X)$ such that $\sum_i\pi_i=\Delta_X$ in
$\Cor_S^0(X,X)$ and the following holds:

\medskip\noindent
(3.3.1)\, The restriction of $(\pi_i)_*$ to
$({}^pR^j(f_{\ok})_*\bQ_l[n])[-j]$ under an isomorphism (3.2.1)
is the inclusion $({}^pR^j(f_{\ok})_*\bQ_l[n])[-j]\to
\bR(f_{\ok})_*\bQ_l[n]$ if $j=i$ and vanishes otherwise.

\medskip
Similarly, we say that $f:X\to S$ admits a {\it relative
Chow-K\"unneth decomposition in the strong sense} (see [4]),
if there are mutually orthogonal idempotents
$\pi_{i,Z}\in\Cor^0_S(X,X)$ such that $\sum_{i,Z}\pi_{i,Z}
=\Delta_X$ in $\Cor_S^0(X,X)$ and the following holds:

\medskip\noindent
(3.3.2)\, The restriction of $(\pi_{i,Z'})_*$ to
$\IC_ZE^j_{Z^{\circ}}[-j]$ under an isomorphism (3.2.3) is the
inclusion $\IC_ZE^j_{Z^{\circ}}[-j]\to\bR(f_{\ok})_*\bQ_l[n]$
if $(j,Z)=(i,Z')$ and vanishes otherwise.

\medskip\noindent
Here we assume $\pi_i=0$ (or $\pi_{i,Z}=0)$ if
${}^pR^i(f_{\ok})_*\bQ_l[n]=0$ (or $\IC_ZE^i_{Z^{\circ}}=0$).
Note that the above conditions (3.3.1--2) may depend on the
choice of the isomorphism (3.2.1) or (3.2.3), and we say that
$f$ admits a relative Chow-K\"unneth decomposition if there is
an isomorphism (3.2.1) or (3.2.3) such that the condition
(3.3.1) or (3.3.2) is satisfied.
By Proposition below, however, these conditions can be
replaced with weaker ones.

\medskip\noindent
{\bf 3.4.~Proposition.} {\it
The conditions {\rm (3.3.1)} and {\rm (3.3.2)} in the above
definition may be replaced respectively by the following\,{\rm :}

\medskip\noindent
{\rm (3.4.1)}\, The induced morphism ${}^p\cH^j(\pi_i)_*\in
\End({}^pR^j(f_{\ok})_*\bQ_l[n])$ is the identity on
${}^pR^j(f_{\ok})_*\bQ_l[n]$ if $i=j$, and vanishes otherwise.

\medskip\noindent
{\rm (3.4.2)}\, The restriction of ${}^p\cH^j(\pi_{i,Z'})_*$ to
$\IC_ZE^j_{Z^{\circ}}\subset {}^pR^j(f_{\ok})_*\bQ_l[n]$
is the inclusion $\IC_ZE^j_{Z^{\circ}}\to
{}^pR^j(f_{\ok})_*\bQ_l[n]$ if $(j,Z)=(i,Z')$,
and vanishes otherwise.
}

\medskip\noindent
{\it Proof.}
For the assertion on $(\pi_i)_*$, we have by Lemma~(3.5)
below a factorization
$$(\pi_i)_*:\bR(f_{\ok})_*\bQ_l[n]\buildrel{\rho_i}\over\to
({}^pR^i(f_{\ok})_*\bQ_l[n])[-i]\buildrel{\iota_i}\over\to
\bR(f_{\ok})_*\bQ_l[n].\leqno(3.4.3)$$
(Note that (3.2.1) obtained in (3.5) may depend on $i$.)
We then get $\rho_i\scirc\iota_j=0$ for $i\ne j$ using
$\pi_i\scirc\pi_j=0$.
This gives an isomorphism (3.2.1) such that the condition (3.3.1)
is satisfied.
The argument is similar for $(\pi_{i,Z})_*$.

\medskip\noindent
{\bf 3.5.~Lemma.} {\it
Let $\alpha$ be an idempotent of $\End(\bR(f_{\ok})_*\bQ_l[n])$
and $i\in\bZ$ such that ${}^p\cH^j\alpha=0$ for $j\ne i$
in $\End({}^pR^j(f_{\ok})_*\bQ_l[n])$.
Let $\cF=\rIm\,{}^p\cH^i\alpha\subset{}^pR^i(f_{\ok})_*\bQ_l[n]$.
Then there is a factorization
$$\alpha:\bR(f_{\ok})_*\bQ_l[n]\buildrel{\rho}\over\to\cF[-i]
\buildrel{\iota}\over\to\bR(f_{\ok})_*\bQ_l[n].\leqno(3.5.2)$$
}

\medskip\noindent
{\it Proof.}
Choose any isomorphism (3.2.1) so that $\alpha$ is expressed by
a matrix
$$A=(a_{r,s})\quad\text{with}\quad a_{r,s}\in\Ext^{s-r}
({}^pR^s(f_{\ok})_*\bQ_l[n],{}^pR^r(f_{\ok})_*\bQ_l[n]).$$
We have $A^2=A$, $a_{r,s}=0$ for $r>s$, and $a_{r,r}=0$ for
$r\ne i$.
These conditions imply $a_{r,s}=0$ if $r,s>i$ or $r,s<i$,
because they define submatrices which are nilpotent and also
idempotent.
Similarly, we see that $a_{i,s}$ and $a_{r,i}$ factor through
$\cF$ (using a refinement of the matrix associated with the
decomposition of ${}^pR^i(f_{\ok})_*\bQ_l[n]$ defined by
the idempotent ${}^p\cH^i\alpha$).

If we change the isomorphism (3.2.1), this corresponds to an
invertible matrix $B=(b_{r,s})$ with $b_{r,r}=id$, where
$A$ is replaced by $B^{-1}AB$.
Using $B=(b_{r,s})$ with $b_{i,s}=-a_{i,s}$ for $s>i$ and
$b_{r,s}= 0$ for $s>r\ne i$, we may assume $a_{i,s}=0$ for
$s>i$, and then $a_{r,i}=0$ for $r<i$ by a similar argument.
Then we get $a_{r,s}=0$ for $r\ne s$ using $A^2=A$.
So (3.5.2) follows.

\bigskip\bigskip
\centerline{\bf 4. Case of relative dimension 1}

\bigskip\noindent
{\bf 4.1.}
Let $f:X\to S$ be a surjective morphism of irreducible smooth
projective varieties over a perfect field $k$ with $\dim X = 2$,
$\dim S = 1$.
Assume that the generic fiber of $f$ is smooth over $k(S)$, and
each closed fiber of $f$ is a divisor with simple normal
crossings.
Let $U$ be a non-empty open subvariety of $S$ over which $f$ is
smooth. Set $X'=f^{-1}(U)$ and $\Sigma=S\setminus U$.

For a closed point $s$ of $S$, let $X_s$ be the fiber of $f$ over
$s$, and $X_{s,\alpha}$ be the connected components of $X_s$.
For $s\in\Sigma$, let $D_{s,i}$ be the irreducible components
of $X_s$ with reduced structure for $i\in I(s)$.
Then $I(s)$ is the disjoint union of $I(s,\alpha)$ such that
$(X_{s,\alpha,})_{\red}=\bigcup_{i\in I(s,\alpha)}D_{s,i}$ for
any $\alpha$.

Calculating the stalks of the direct image,
we have in the notation of (3.2)
$$\IC_ZE^j_{Z^{\circ}}=0\quad\text{unless}\,\,Z=S,|j|\le 1\,\,
\text{or}\,\,Z=\{s\}\subset\Sigma,j=0.\leqno(4.1.1)$$
Here we may put $Z^{\circ}=U$ if $Z=S$.
Note that $Z^{\circ}=\{s\}$ if $Z=\{s\}\subset\Sigma$.

Let $\tS=\cSpec_S\,f_*\cO_X$ so that we have the Stein
factorization
$$f:X\buildrel{\tf}\over\to\tS\to S.\leqno(4.1.2)$$
Note that $\tS$ is $1$-dimensional and normal, and is hence
regular and smooth over $k$ because $k$ is a perfect field.

\medskip\noindent
{\bf 4.2.~Remarks.}
(i) With the notation and the assumptions of (4.1), let
$[D_{s,i}]$ denote the cycle class of $D_{s,i}$ in
$\CH^2(X)_{\bQ}$.
Let $[X_{s,\alpha}]=\sum_{i\in I(s,\alpha)}m_{s,i}[D_{s,i}]$,
where $m_{s,i}$ is the multiplicity of $X_s$ at a generic
point of $D_{s,i}$.
Then it is well known that the intersection form on
$V_{\alpha}=:\mopl_{i\in I(s,\alpha)}\bQ[D_{s,i}]$ is negative
semi-definite with $[X_{s,\alpha}]\cdot\xi=0$
for any $\xi\in V_{\alpha}$ and the induced bilinear form on
$V'_{\alpha}:=V_{\alpha}/\bQ[X_{s,\alpha}]$ is negative
definite.

Indeed, using the Stein factorization and replacing $S$ if
necessary, we may assume that the $X_s$ are connected.
Then $[X_s]\cdot\xi=0$ moving the 0-cycle $\{s\}$ on $S$.
Let $v_i=[D_{s,i}]$, $w=[X_s]$, and $a_{i,j}=v_i\cdot v_j$.
Replacing $v_i$ if necessary,
we may assume $w=\sum_iv_i$ so that $\sum_i a_{i,j}=0$.
Then we have
$$\msum_i c_iv_i\cdot\msum_j c_j v_j=\msum_{i,j}a_{i,j}c_i c_j
=-\msum_{i\ne j}a_{i,j}(c_i-c_j)^2/2.$$

\medskip
(ii) With the notation and the assumptions of (4.1), assume $k$
is algebraically closed.
Then for $s\in\Sigma$ we have a canonical isomorphism
$$E^0_{\{s\}}=\mopl_{\alpha}\bigl((\mopl_{i\in I(s,\alpha)}
\bQ_l[D_{s,i}])/\bQ_l[X_{s,\alpha}]\bigr),\leqno(4.2.1)$$
where $[D_{s,i}]$ is the cycle class of $D_{s,i}$ in
$H^2(X,\bQ_l)$,
and $[X_{s,\alpha}]=\sum_{i\in I(s,\alpha)}m_{s,i}[D_{s,i}]$
as in Remark (i) above.

Indeed, using the Stein factorization, the assertion can be
reduced to the case where the $X_s$ are connected so that
$\IC_SE^{j}_U$ for $|j|=1$ is a constant perverse sheaf
$\bQ_l[1]$ on $S$ up to a Tate twist,
because the intermediate direct images commute with the
direct image by a finite morphism.
Then we get the assertion using the local cohomology
$H^2_{X_s}(X,\bQ_l)$ (which is a free free $\bQ_l$-module with
basis $[D_{s,i}]$), because the image of
$[X_s]$ in $H^2(X,\bQ_l)$ coincides with that of $[X_{s'}]$ for
$s'\in U$, and hence with $H^2(S,\IC_SE^{-1}_U)$.
Note that the latter is naturally identified with a subspace of
$H^2(X,\bQ_l)$ by the $E_2$-degeneration of the (perverse) Leray
spectral sequence associated with the perverse $t$-structure,
see [2].

\medskip
The following is a variant of a result in [9] where
$X,S$ may be higher dimensional but the singular fiber are
more restrictive.

\medskip\noindent
{\bf 4.3.~Proposition.} {\it
With the notation and the assumptions of {\rm (4.1)}, 
let $\pi'\in\Cor^0_U(X',X')$ be an idempotent whose action
on ${}^pR^jf_*\bQ_l[2]|_U$ vanishes for $j\ne 0$.
Then $\pi'$ can be extended uniquely to an idempotent $\pi$ of
$\Cor^0_S(X,X)$ such that the action of ${}^p\cH^j\pi_*$
vanishes on ${}^pR^j(f_{\ok})_*\bQ_l[2]$ for $j\ne 0$ and
also on $\IC_{\{s\}}E^0_{\{s\}}$ for $s\in\Sigma$ and $j=0$.
}

\medskip\noindent
{\it Proof.}
Let $\pi$ be any extension of $\pi'$ to $\Cor^0_S(X,X)=
\CH^1(X\mtim_SX)_{\bQ}$.
Since the ambiguity of $\pi$ is given by an element of
$\mopl_{s\in\Sigma}\CH^0(X_s\mtim_{k(s)}X_s)_{\bQ}$,
we may assume that the restriction of ${}^p\cH^0\pi_*$
to $\IC_{\{s\}}E^0_{\{s\}}$ vanishes by (4.2).
The action of ${}^p\cH^j\pi_*$ on ${}^pR^j(f_{\ok})_*\bQ_l[2]$
vanishes for $|j|=1$ because these sheaves are constant
perverse sheaves $\bQ_l[1]$ on $S$ up to a Tate twist.

Set $\gamma=\pi^2-\pi$. Since $\pi'{}^2=\pi'$, we have
$\gamma\in\mopl_{s\in\Sigma}\CH^0(X_s\mtim_{k(s)}X_s)_{\bQ}$,
and the restriction of ${}^p\cH^0\gamma_*$ to
$\IC_{\{s\}}E^0_{\{s\}}$ vanishes by the above condition on
${}^p\cH^0\pi_*$.
Combined with (4.2) this implies that
$$\gamma=\msum_{s\in\Sigma}\msum_{\alpha}\bigl([X_{s,\alpha}]
\mtim\xi_{s,\alpha}+\xi'_{s,\alpha}\mtim[X_{s,\alpha}]\bigr)
\quad\text{with}\quad\xi_{s,\alpha},\xi'_{s,\alpha}\in
\CH^0(X_s)_{\bQ}.$$

For $\xi\in\CH^0(X_s)_{\bQ}$ with $s\in\Sigma$, we have
$$\bigl([D_{s,i}]\mtim\xi\bigr)\sbull\pi=[D_{s,i}]
\mtim\xi'\quad\text{with}\quad\xi'\in\CH^0(X_s)_{\bQ},$$
and the image of $\xi'$ in $\CH^1(X)_{\bQ}$ coincides with
the image of $\xi$ by the action of $\pi\in\Cor_k^0(X,X)$
(see (3.1.1) for the definition of $\sbull$).
So we get $\xi'\in\msum_{\beta}\bQ[X_{s,\beta}]$ by the
hypothesis on the action of $\pi_*$ together with (4.2).
(Note that if $\eta\in\CH^1(X)_{\bQ}$ is numerically equivalent
to nonzero, then the map $\CH^1(X)_{\bQ}\to\CH^2(X\mtim_kX)
_{\bQ}$ defined by $\xi'\mapsto\eta\mtim\xi'$ is injective
using a correspondence defined by $\eta'\in\CH^1(X)_{\bQ}$ such
that $\eta\cdot\eta'\ne 0$.)

For $\pi\sbull\bigl([X_{s,\alpha}]\mtim\xi\bigr)$, we see that
it coincides with $\eta\mtim\xi$ with $\eta$ the restriction
to $X_s$ of $(pr_1)_*\pi\in\CH^{0}(X\mtim_S\tS)$,
where $pr_1:X\mtim_SX\to X\mtim_S\tS$ is the base change of
$\tf:X\to\tS$, and the inclusion $X_s\to X\mtim_S\tS$ is the
base change of $s'\to\tS$ corresponding to $X_{s,\alpha}$.
Then the pushdown of cycle $(pr_1)_*\pi$ vanishes in
$\CH^{0}(X\mtim_S\tS)$ because ${}^p\cH^{-1}\pi_*=0$.
(Here we can restrict these over $U$.)
So we get $\pi\sbull\bigl([X_{s,\alpha}]\mtim\xi\bigr)=0$.

We have similar assertions for $\xi\mtim[X_{s,\alpha}]$.
So, replacing $\pi$ with $\pi+\gamma$, the assertion is reduced
to the case $\gamma\in\msum_{s\in\Sigma}\msum_{\alpha,\beta}
\bQ[X_{s,\alpha}]\mtim[X_{s,\beta}]$.
Then, repeating the same argument and replacing $\pi$ as above,
we get $\gamma=0$.

The proof of the uniqueness is similar.
If $\pi+\gamma$ is also an idempotent whose action satisfies the
above conditions, we first get $\gamma\in\msum_{s\in\Sigma}
\msum_{\alpha,\beta}\bQ[X_{s,\alpha}]\mtim[X_{s,\beta}]$,
and then $\gamma=0$.

\medskip\noindent
{\bf 4.4.~Proposition.} {\it
With the notation and the assumptions of {\rm (4.1)}, 
let $\pi\in\Cor^0_S(X,X)$ such that the action of
${}^p\cH^j\pi_*$ vanishes on ${}^pR^j(f_{\ok})_*\bQ_l[2]$ for
$j\ne 0$ and also on $\IC_{\{s\}}E^0_{\{s\}}$ for $s\in\Sigma$
and $j=0$.
Let $\gamma\in\Cor^0_S(X,X)$ which is represented by a cycle
supported on $X_s\mtim_{k(s)}X_s$ for a closed point $s$ of
$S$. Then $\pi^2\scirc\gamma=\gamma\scirc\pi^2=0$.
}

\medskip\noindent
{\it Proof.}
This follows from the same argument as in the proof of
Proposition~(4.3).

\medskip
The following gives an absolute decomposition of a
relative decomposition in the surface case.

\medskip\noindent
{\bf 4.5.~Proposition.} {\it
With the notation and the assumptions of {\rm (4.1)}, 
let $\pi\in\Cor^0_S(X,X)$ be an idempotent such that the
action of ${}^p\cH^j\pi_*$ on ${}^pR^j(f_{\ok})_*\bQ_l[2]$
vanishes for $j\ne 0$.
Then there are mutually orthogonal idempotents
$\pi_i\in\Cor^0_k(X,X)$ for $i=-1,0,1$ such that
$\sum_i\pi_i=\pi$ in $\Cor^0_k(X,X)$ and the action of
$H^{j+2}(\pi_i)_*$ on $H^{j+2}(X_{\ok},\bQ_l)$ vanishes
for $j\ne i$.
}

\medskip\noindent
{\it Proof.}
Let $s$ be a general closed point of $S$, and put $Y=X_s,
k'=k(s)$.
We denote the restriction of $\pi$ to $Y\mtim_{k'}Y$ by
$$\pi'\in\Cor_{k'}^0(Y,Y)=\CH^1(Y\mtim_{k'}Y)_{\bQ}.$$
Its image in $\Cor_{k}^0(Y,Y)$ will be denoted also
by $\pi'$.
The action of $H^j\pi'_*$ on $H^j(Y_{\ok},\bQ_l)$ vanishes
for $j\ne 1$ by the assumption on the action of
${}^p\cH^{j-1}\pi_*$ together with the generic base change
theorem [7], where $Y_{\ok}=Y\otimes_k{\ok}$ (and not
$Y\otimes_{k'}{\ok})$.

By an argument similar to the proof of Proposition~(2.3),
we have the projector $\pi''\in\Cor_{k}^0(Y,Y)$ such
that the action of $H^j\pi''_*$ on $H^j(Y_{\ok},\bQ_l)$
vanishes for $j\ne 1$ and the image of $H^1\pi''_*$
coincides with that of the restriction morphism
$H^1(X_{\ok},\bQ_l)\to H^1(Y_{\ok},\bQ_l)$.
The last morphism is also induced by $\Spec\,k(s)\to S$, and
is compatible with the action of $\pi$. This holds also for
the Gysin morphism. So the actions of $H^1\pi'_*$ and
$H^1\pi''_*$ on $H^1(Y_{\ok},\bQ_l)$ commute, and we get
a projector $\pi'\scirc\pi''\scirc\pi'\in\Cor_{k}^0(Y,Y)$.
This induces a projector of restricted type $\pi'_{-1}\in
\Cor_k^0(X,X)$ by the same argument as in (2.4). Then
$\pi_{-1}:=\pi\scirc\pi'_{-1}\scirc\pi$ is also a projector,
and is a refinement of $\pi$.
The argument is similar for $\pi_1$. We can modify it
so that $\pi_1$ is orthogonal to $\pi_{-1}$ by (1.4).

\bigskip\bigskip
\centerline{\bf 5. Construction of relative projectors}

\bigskip\noindent
{\bf 5.1.~Proof of Theorem~2.}
Since $f$ has only isolated singularities, we have in the
notation of (3.2)
$$\IC_ZE^j_{Z^{\circ}}=0\quad\text{unless}\,\,Z=\{s\}\subset
\Sigma\,\,\text{with}\,\,j=0,\,\,\text{or}\,\,Z=S.
\leqno(5.1.1)$$
This follows form the commutativity of the vanishing cycle
functor $\varphi$ (see [6]) with the direct image under a proper
morphism, because the assumption implies that
$\varphi_{f^*t}\bQ_{l}[n]$ is a perverse sheaf with a
$0$-dimensional support where $t$ is a local coordinate on $S$.

Assume $j=2p+1<m$. By (2.1-3), there is a smooth projective
variety $C$ of pure dimension $1$ over $K$ with a correspondence
$\Gamma\in\CH^{p+1}(C\mtim_KY)_{\bQ}$ inducing the surjection
$$\Gamma_*:H^1(C_{\oK},\bQ_l)(-p)\to H^j(Y_{\oK},\bQ_l),$$
where $C_{\oK}=C\motim_K\oK$.
Let $Z$ be a smooth projective variety of pure dimension $2$ over
$k$ with a flat morphism $g:Z\to S$ whose generic fiber is
isomorphic to $C$. This is obtained by using a compactification of
a spreading-out of $C$ together with a resolution of singularities
in the 2-dimensional case.
Here we may assume that each closed fiber of $g$ is a divisor
with normal crossings.

Taking the closure of a spreading-out of $\Gamma$, we get
$\Gamma_S\in\CH^{p+1}(Z\mtim_SX)_{\bQ}$ whose base change
by $\Spec\,K\to S$ is $\Gamma$.
Then we get the relative projector $\pi_{j-m,S}$
using arguments similar to (2.4) together with
(4.3). We can modify these so that they are orthogonal
to each other using (4.4).
The argument is similar for the case $j=2p$.
Note that $\pi_0$ is defined by $\Delta_X-\msum_{j\ne 0}\pi_j$.

To get the absolute Chow-K\"unneth decomposition, we
have a refinement of the relative projector $\pi_{j,S}$
using (4.5).

\end{document}